\documentclass{amsart}

\usepackage[T1]{fontenc}
\usepackage[utf8]{inputenc}
\usepackage{amssymb,mathtools}
\usepackage[hidelinks]{hyperref}
\hypersetup{pdfauthor={Damián Pinasco, Juan Pablo Pinasco, María Victoria Venuti},pdftitle={Optimal polynomial meshes on quasicircles via conformal scale comparison},pdfsubject={},pdfkeywords={}}

\newtheorem{theorem}{Theorem}[section]
\newtheorem{proposition}[theorem]{Proposition}
\newtheorem{lemma}[theorem]{Lemma}
\newtheorem{corollary}[theorem]{Corollary}
\theoremstyle{definition}
\newtheorem{definition}[theorem]{Definition}
\theoremstyle{remark}
\newtheorem{remark}[theorem]{Remark}

\newcommand{\C}{\mathbb C}
\newcommand{\R}{\mathbb R}
\newcommand{\T}{\mathbb T}
\newcommand{\dist}{\operatorname{dist}}
\newcommand{\diam}{\operatorname{diam}}
\newcommand{\card}{\operatorname{card}}
\newcommand{\norm}[1]{\left\lVert #1\right\rVert}

\begin{document}

\title[Optimal polynomial meshes on quasicircles]{Optimal polynomial meshes
on quasicircles via conformal scale comparison}

\author{Damián Pinasco}
\address{Departamento de Matemática y Estadística, Universidad Torcuato Di Tella, Av. Figueroa Alcorta 7350 (C1428BCW), Ciudad de Buenos Aires, Argentina; CONICET, Buenos Aires, Argentina}
\email{dpinasco@utdt.edu}

\author{Juan Pablo Pinasco}
\address{IMAS (UBA-CONICET), Departamento de Matemática, Facultad de Ciencias Exactas y Naturales, Universidad de Buenos Aires, Intendente Güiraldes 2160, Ciudad Universitaria, Pabellón I, C1428EGA Buenos Aires, Argentina}
\email{jpinasco@dm.uba.ar}

\author{María Victoria Venuti}
\address{Departamento de Matemática y Estadística, Universidad Torcuato Di Tella, Av. Figueroa Alcorta 7350 (C1428BCW), Ciudad de Buenos Aires, Argentina}
\email{victoria.venuti@utdt.edu}

\thanks{Accepted for publication in the Proceedings of the American Mathematical Society.}

\subjclass[2020]{Primary 41A17; Secondary 30C62, 30C85, 30E10, 41A10}
\keywords{Admissible polynomial mesh, sampling discretization, norming set,
quasicircle, quasiconformal map, harmonic measure, von Koch snowflake}

\begin{abstract}
We prove that every quasicircle $\Gamma\subset\C$ admits an optimal admissible polynomial
mesh for univariate complex polynomials. More precisely, if $E$ denotes the closed Jordan domain
bounded by $\Gamma$, then for every $\Lambda>1$ there are finite sets $A_n\subset\Gamma$
with
\[
  \card(A_n)=O(n),
  \qquad
  \norm p_E=\norm p_\Gamma\le \Lambda\norm p_{A_n},
  \qquad p\in\mathcal P_n(\C).
\]
The nodes may be chosen as images of equally spaced points of the unit circle under the
inverse exterior conformal map, equivalently as equal-mass points for harmonic measure at infinity.
The proof combines the classical conformal level-curve scale with a discrete norming
principle. If $\rho_n(z)$ denotes the distance from a boundary point to the exterior level
curve $|\Phi|=1+1/n$, Bernstein--Walsh and Cauchy's formula show that any sufficiently fine
$\rho_n$-net is norming. A scale-comparison property then turns $O(n)$ equally spaced
angular samples into such a net. Quasicircles satisfy this property by quasisymmetry. The
argument requires no smoothness or rectifiability and applies in particular to the
classical von Koch snowflake.
\end{abstract}

\maketitle

\section{Introduction}

\subsection{Polynomial meshes in one complex variable}

Let $E\subset\C$ be compact and
\[
  \mathcal P_n(\C)=\{p(z)=a_0+a_1z+\cdots+a_nz^n\}.
\]
A sequence of finite sets $A_n\subset E$ is an \emph{admissible polynomial mesh} (AM) if
there is a constant $\Lambda\ge1$, independent of $n$, such that
\begin{equation}\label{eq:introAM}
  \norm p_E\le \Lambda\norm p_{A_n},
  \qquad p\in\mathcal P_n(\C),
\end{equation}
and $\card(A_n)=O(n^s)$ for some $s>0$. The terminology and the modern mesh framework
go back to Calvi and Levenberg \cite{CalviLevenberg}. If $E$ is
polynomially determining, then
\[
  \dim_\C\bigl(\mathcal P_n(\C)|_E\bigr)=n+1,
\]
so every norming set for $\mathcal P_n(\C)|_E$ contains at least
$n+1$ points. In this setting, an AM with
\[
  \card(A_n)=O(n)
\]
is therefore called \emph{optimal}.

The theory of optimal polynomial meshes has also been developed in the real
multivariate setting, where optimal meshes have $O(n^d)$ points; see, among
others, Kro\'o \cite{Kroo2011}, Piazzon \cite{Piazzon2016} and Dai--Prymak
\cite{DaiPrymak}. In particular, Dai and Prymak proved that every convex body
in $\R^d$ admits an optimal mesh. Here we consider instead polynomials of one
complex variable on planar compact sets, for which the dimension of the
degree-$n$ polynomial space is $n+1$.

A general near-optimal result is known with essentially no geometric regularity. Bloom,
Bos, Calvi and Levenberg proved that every polynomially determining compact subset of
$\C^d$ has an AM of cardinality $O((n\log n)^d)$
\cite[Proposition~23]{BloomBosCalviLevenberg}. In one complex variable, every infinite
compact is polynomially determining, hence every infinite $E\subset\C$ admits an AM with
$O(n\log n)$ points.

Removing the logarithm is a geometric problem. Related constructions
for regular curves use conformal or analytic transformations and
perturbation arguments. Bia\l as-Cie\.z and Calvi
\cite[Section~4.2]{BialasCiezCalvi} obtained weakly admissible
meshes on compact sets bounded by Alper-smooth Jordan curves from
conformal images of equally spaced points of the unit circle, with
$d+1$ nodes in degree $d$ and mesh constants of order $O(\log d)$. Analytic transformations
and perturbations of polynomial meshes were studied by
Piazzon and Vianello
\cite{PiazzonVianello2010,PiazzonVianello}. In particular, the
perturbation method of \cite{PiazzonVianello} gives meshes of cardinality
$O(n^{1+2/k})$ under $C^k$ regularity, $k>2$. These are outputs of particular methods and
are not sharpness statements for all Jordan boundaries. For a $C^2$ Jordan curve, the Bernstein inequality has the
linear degree scale
\[
  \norm{p'}_\Gamma\le C_\Gamma n\norm p_\Gamma;
\]
see, for instance, \cite{KalmykovNagyTotik}. Consequently,
arclength sampling at step $O(1/n)$ yields an optimal mesh. Explicit optimal meshes for boundaries formed by a fixed number
of polynomial or trigonometric polynomial arcs were constructed by Bia\l as-Cie\.z, Kenne,
Sommariva and Vianello \cite{BialasKenneSommarivaVianello}. Such constructions include
polygons and many curvilinear or cuspidal examples.

Thus linear meshes are known in several regular or explicitly parametrized situations,
while $O(n\log n)$ is available abstractly for every infinite compact. The gap between
these two regimes leaves open the role of geometry in forcing linear cardinality, even
in one complex variable.

A different recent approach by Pinasco and Venuti
\cite{PinascoVenutiGreen} uses Green sublevels to construct
optimal polynomial meshes for classes of compact subsets of
the real line, including highly disconnected examples. Both approaches use Bernstein--Walsh to convert control on a
nearby Green level into a norming estimate; here the counting
mechanism is instead provided by the conformal angular
parameter on a Jordan boundary.

\subsection{Constructive approximation on quasidisks}

There is a separate and much older literature that is directly relevant to the geometry of
the present problem. A \emph{quasidisk} is a Jordan domain whose boundary is a quasicircle.
Polynomial approximation on quasidisks and quasiconformal arcs has been studied since the
work of Dzjadyk, Belyi and Andrievskii and their collaborators. Belyi used conformal and
quasiconformal methods to obtain approximation theorems on domains with quasiconformal
boundary \cite{Belyi1977}; Andrievskii extended direct theorems to arbitrary
quasiconformal arcs, including everywhere nonrectifiable ones \cite{Andrievskii1981}.

A central object in this theory is the distance from a boundary point to a nearby exterior
conformal level curve. In the notation used below, this is precisely the scale
\[
  \rho_n(z)=\dist\bigl(z,\{ |\Phi|=1+1/n\}\bigr).
\]
Its role is classical. In particular, Andrievskii formulated direct approximation theorems
in terms of this distance \cite{Andrievskii1986}; see also the monograph
\cite{AndrievskiiBelyiDzjadyk} and the survey \cite{AndrievskiiSurvey2002}. The same
circle of ideas underlies later work on simultaneous approximation and interpolation
\cite{AndrievskiiPritskerVarga}, maximal polynomial ranges \cite{AndrievskiiMaxRange},
maximum-modulus refinements on quasidisks \cite{AndrievskiiMaxModulus}, and the
Totik--Widom property for quasidisks \cite{AndrievskiiNazarov}.

The same scale also appears in distortion estimates
for exterior conformal maps of quasidisks and in
Bernstein--Dzjadyk inequalities of the form
\[
  |p'(z)|\le C_\Gamma
  \frac{\norm p_E}{\rho_n(z)},
  \qquad z\in\Gamma,\quad p\in\mathcal P_n(\C);
\]
see \cite{AndrievskiiBelyiDzjadyk,AndrievskiiSurvey2002}.
Our use of the scale is different: Proposition~\ref{prop:local}
gives a self-contained difference estimate, and
Corollary~\ref{cor:rhonet} turns it into a finite norming criterion.

These classical conformal and quasiconformal tools provide the geometric background for the present work. The question addressed here is whether this structure can be turned into a discrete norming mechanism with optimal cardinality.

More precisely, we ask whether the boundary of a quasidisk can be replaced by only $O(n)$ sample points while preserving a norming inequality uniformly over all polynomials of degree at most $n$. The natural candidates are the boundary points obtained by transporting equally spaced points of the unit circle under the inverse exterior conformal map, namely the classical Fej\'er points introduced in this setting by Fej\'er \cite{Fejer1918} and subsequently studied in complex interpolation, for example by Curtiss \cite{Curtiss1960}.
Koosis \cite{Koosis2001} used nodes equally spaced in exterior harmonic measure on sufficiently smooth Jordan curves to estimate polynomial values at interior points through logarithmic sums.

Related work of Chui and Zhong \cite{ChuiZhong2021}
studies Marcinkiewicz--Zygmund inequalities on $L$-shape
arcs using the same exterior level-curve scale.
Their results concern $L^p$ inequalities, $1<p<\infty$,
whereas the present problem is a uniform-norm norming
inequality with oversampling; in the terminology of sampling
discretization \cite{KKLT2022}, it asks for a discretization of the
uniform norm using a number of sampling points proportional to the
dimension.

We prove that suitably oversampled Fej\'er points yield
admissible meshes of optimal linear cardinality on every quasicircle. This extends the linear-mesh regime from the regular and explicitly parametrized settings treated above to a class that includes nonrectifiable and fractal boundaries. The key step is to combine the local conformal scale with a scale-comparison principle that turns equally spaced angular samples into an admissible polynomial mesh of cardinality $O(n)$.

\subsection{The scale-comparison mechanism}

The argument has two parts. First, Bernstein--Walsh and Cauchy's formula give a
geometry-free local principle: any boundary set which is sufficiently dense relative to
$\rho_n(z)$ is norming for $\mathcal P_n(\C)$. Second, we isolate a scale-comparison
property (SC) which allows such a variable-radius net to be constructed with $O(n)$ points.
For quasicircles, (SC) follows from one application of quasisymmetry. Thus the logical
structure is
\[
  \text{quasicircle}
  \Longrightarrow \text{(SC)}
  \Longrightarrow \rho_n\text{-net with }O(n)\text{ points}
  \Longrightarrow \text{optimal mesh}.
\]
The middle condition is stated entirely in terms of the exterior conformal parameter; the
present paper uses it as a sufficient hypothesis for the sampling argument.

This formulation also makes clear that the result is not perturbative. A quasicircle need
not be close to a circle in Hausdorff distance, arclength, Lipschitz geometry, or any smooth
sense; it may be nonrectifiable and fractal. What is controlled is relative scale distortion.
Equal angular increments can have very different Euclidean sizes, but an angular increment
of size $c/n$ is uniformly small compared with the local exterior scale $\rho_n(z)$ when
$c$ is small.

\subsection{Organization of the paper}

Section~\ref{sec:prelim} recalls the complex-analytic and quasiconformal background and
places the conformal level-curve scale in its classical setting. Section~\ref{sec:rho}
proves the local Green-scale norming principle. Section~\ref{sec:SC} introduces (SC) and
shows that it produces optimal meshes. Section~\ref{sec:QC} proves that every quasicircle
satisfies (SC). Section~\ref{sec:koch} records the von Koch consequence, and the final
section contains some concluding remarks.

\section{Preliminaries}\label{sec:prelim}

\subsection{Jordan domains, Green functions and Bernstein--Walsh}

From now on, let $E$ be a closed Jordan domain and write
$\Gamma:=\partial E$. Put
\[
  \Omega:=\widehat\C\setminus E.
\]
The normalized exterior conformal map
\[
  \Phi:\Omega\longrightarrow\widehat{\C}\setminus\overline{\mathbb D},
  \qquad
  \Phi(\infty)=\infty,
\qquad
\lim_{z\to\infty}\frac{\Phi(z)}{z}>0,
\]
is uniquely determined; write $\Psi:=\Phi^{-1}$.
By Carath\'eodory's theorem, $\Psi$ extends continuously to $\T$,
and its boundary restriction is a homeomorphism from $\T$ onto $\Gamma$;
see, for example, Pommerenke \cite{Pommerenke}.

The Green function of $\Omega$ with pole at infinity is
\[
  g_E(z)=\log|\Phi(z)|,
  \qquad z\in\Omega.
\]
Since $|\Phi(z)|\to1$ as $z\to\Gamma$, this function extends continuously
to $\C$ by setting $g_E=0$ on $E$. For $p\in\mathcal P_n(\C)$,
Bernstein--Walsh gives
\begin{equation}\label{eq:BW}
  |p(z)|\le \norm p_E e^{n g_E(z)},
  \qquad z\in\C,
\end{equation}
with the estimate on $E$ understood trivially. We refer to Ransford \cite{Ransford} for the
potential-theoretic background. The maximum principle yields
\begin{equation}\label{eq:maxprinciple}
  \norm p_E=\norm p_\Gamma.
\end{equation}

The use of distances from exterior conformal level curves is classical in constructive
complex approximation. In particular, the scale at level $1+1/n$ already plays a central
role in the work of Belyi and Andrievskii \cite{Belyi1977,Andrievskii1986}; see also
\cite{AndrievskiiBelyiDzjadyk,AndrievskiiSurvey2002}. We use that classical scale for a
different purpose: to formulate a finite norming criterion.

\subsection{Quasicircles and quasisymmetry}

A Jordan curve $\Gamma$ is a \emph{$K$-quasicircle} if $\Gamma=F(\T)$ for a
$K$-quasiconformal homeomorphism $F$ of the Riemann sphere. We shall use two standard
consequences; see Ahlfors \cite{AhlforsReflection}, Gehring--Hag \cite{GehringHag}, and
Heinonen \cite{Heinonen}.

First, quasicircles are exactly the Jordan curves of bounded turning: there is $B<\infty$
such that, for every $a,b\in\Gamma$, the smaller-diameter subarc $\Gamma[a,b]$ satisfies
\begin{equation}\label{eq:BT}
  \diam\Gamma[a,b]\le B|a-b|.
\end{equation}
Second, the inverse exterior map admits a quantitative
quasiconformal extension. Indeed, let
\[
  J(w):=\frac{1}{\overline w}.
\]
If $\Gamma=F(\T)$ with $F$ $K$-quasiconformal, then
\[
  R:=F\circ J\circ F^{-1}
\]
is a $K^2$-quasiconformal reflection across $\Gamma$.
Hence the standard extension by reflection
\[
  \widetilde\Psi(w):=
  \begin{cases}
    \Psi(w), & |w|\ge1,\\[2mm]
    R\bigl(\Psi(J(w))\bigr), & |w|<1,
  \end{cases}
\]
defines a global $K^2$-quasiconformal homeomorphism of
$\widehat\C$ which agrees with $\Psi$ on $\{|w|\ge1\}$. Since $\widetilde\Psi(\infty)=\infty$, its restriction to $\C$
is a $K^2$-quasiconformal homeomorphism of the plane. Hence,
by \cite[Theorem~11.14]{Heinonen}, applied with domain and target
the whole plane, it is
$\eta_K$-quasisymmetric, where $\eta_K$ depends only on $K$.
Thus there exists an increasing homeomorphism
\[
  \eta_K:[0,\infty)\longrightarrow[0,\infty),
  \qquad \eta_K(0)=0,
\]
such that
\begin{equation}\label{eq:QS}
  \frac{|x-a|}{|x-b|}\le t
  \quad\Longrightarrow\quad
  \frac{|\widetilde\Psi(x)-\widetilde\Psi(a)|}
       {|\widetilde\Psi(x)-\widetilde\Psi(b)|}
  \le \eta_K(t),
\end{equation}
for all $x,a,b\in\C$ with $x\ne b$. In particular,
\begin{equation}\label{eq:eta0}
  \eta_K(t)\longrightarrow0
  \qquad (t\downarrow0).
\end{equation}
For standard formulations of quasiconformal $\Rightarrow$ quasisymmetric in Euclidean
spaces see also Tukia--V\"ais\"al\"a \cite{TukiaVaisala}.

\section{A local Green-scale norming principle}\label{sec:rho}

For $n\ge1$ set
\begin{equation}\label{eq:rn}
  r_n:=1+\frac1n,
  \qquad
  \Gamma_n:=\Psi(r_n\T),
\end{equation}
and, for $z\in\Gamma$, define
\begin{equation}\label{eq:rhon}
\rho_n(z):=\dist(z,\Gamma_n).
\end{equation}
Since $\Gamma_n$ is compact and disjoint from $\Gamma$,
we have $\rho_n(z)>0$ for every $z\in\Gamma$.

This is the classical local conformal scale associated with the level $1+1/n$; our use of
it below is as a sampling radius rather than as an approximation modulus. Since
$g_E=\log|\Phi|$ on $\Omega$,
\begin{equation}\label{eq:level}
  \Gamma_n=\{z\in\Omega:g_E(z)=\log r_n\}.
\end{equation}
No geometric regularity of $\Gamma$ is used in the following local Green-scale estimate beyond the Jordan
setting above.

\begin{proposition}\label{prop:local}
Let $n\ge1$, $p\in\mathcal P_n(\C)$, $z\in\Gamma$, and $\rho=\rho_n(z)$. If
$|\zeta-z|<\rho$, then
\begin{equation}\label{eq:localdifference}
  |p(\zeta)-p(z)|
  \le
  e\,\norm p_E\,
  \frac{|\zeta-z|}{\rho-|\zeta-z|}.
\end{equation}
\end{proposition}

\begin{proof}
If $|u-z|<\rho$, then $g_E(u)<\log r_n$. Indeed, otherwise continuity of $g_E$ on the
segment from $z$ to $u$ would produce a point of the level set \eqref{eq:level} closer to
$z$ than $\rho$. Hence Bernstein--Walsh gives
\[
  |p(u)|\le r_n^n\norm p_E\le e\norm p_E
\]
throughout $B(z,\rho)$.

Fix $|\zeta-z|<r<\rho$. Cauchy's formula for the difference gives
\[
  p(\zeta)-p(z)
  =\frac{1}{2\pi i}\int_{|u-z|=r}
  p(u)\frac{\zeta-z}{(u-\zeta)(u-z)}\,du.
\]
Therefore
\[
  |p(\zeta)-p(z)|
  \le e\norm p_E\frac{|\zeta-z|}{r-|\zeta-z|}.
\]
Letting $r\uparrow\rho$ proves \eqref{eq:localdifference}.
\end{proof}

The local estimate immediately yields the following $\rho_n$-net principle.

\begin{corollary}\label{cor:rhonet}
Fix $0<q<1$ and set
\begin{equation}\label{eq:tau}
  \tau:=\frac{q}{e+q}.
\end{equation}
Suppose that $A_n\subset\Gamma$ has the following property: for every $z\in\Gamma$ there
exists $a\in A_n$ such that
\begin{equation}\label{eq:rhonet}
  |a-z|\le \tau\rho_n(z).
\end{equation}
Then
\begin{equation}\label{eq:normrho}
  \norm p_E=\norm p_\Gamma
  \le \frac{1}{1-q}\norm p_{A_n},
  \qquad p\in\mathcal P_n(\C).
\end{equation}
\end{corollary}

\begin{proof}
Choose $z_*\in\Gamma$ with $|p(z_*)|=\norm p_\Gamma$ and choose $a\in A_n$ satisfying
\eqref{eq:rhonet}. Proposition~\ref{prop:local} gives
\[
  |p(a)-p(z_*)|
  \le e\norm p_E\frac{\tau}{1-\tau}
  =q\norm p_E.
\]
Using \eqref{eq:maxprinciple},
\[
  |p(a)|\ge (1-q)\norm p_\Gamma,
\]
which is \eqref{eq:normrho}.
\end{proof}

\begin{remark}
Corollary~\ref{cor:rhonet} contains the analytic core of the argument. It does not ask for
a uniform Euclidean mesh. The radius $\rho_n(z)$ is allowed to vary strongly with the base
point. Geometry enters only through the problem of constructing and counting a
$\tau\rho_n$-net.
Related constructions of norming meshes from Bernstein-type inequalities also occur in the real setting; see \cite{Vianello2014} and, for a multivariate construction involving the distance to the boundary, \cite{Piazzon2016}.
\end{remark}

\section{The scale-comparison property}\label{sec:SC}

\begin{definition}\label{def:SC}
A Jordan curve $\Gamma$ is said to satisfy the \emph{scale-comparison property} (SC) if
there exist $c_0>0$ and a function
\[
  \omega_\Gamma:(0,c_0]\longrightarrow(0,\infty),
  \qquad
  \omega_\Gamma(c)\longrightarrow0\quad(c\downarrow0),
\]
such that, for every $n\ge1$, every $\theta\in\R$, every $0<c\le c_0$, and every
$h\in\R$ with $|h|\le c/n$,
\begin{equation}\label{eq:SC}
  |\Psi(e^{i(\theta+h)})-\Psi(e^{i\theta})|
  \le
  \omega_\Gamma(c)\,
  \rho_n(\Psi(e^{i\theta})).
\end{equation}
\end{definition}

Property (SC) is deliberately stated without quasiconformal terminology. It compares two
scales attached to the same conformal coordinate: a tangential angular increment of size
$c/n$ and the Euclidean distance from the boundary to the exterior level curve of conformal
height $1/n$. The following optimal-mesh result shows that this comparison is sufficient for
the mesh argument.

\begin{theorem}\label{thm:SC}
Let $\Gamma$ be a Jordan curve satisfying \textup{(SC)},
and let $E$ be the closed Jordan domain bounded by $\Gamma$.
Then for every $\Lambda>1$ there exist finite sets $A_n\subset\Gamma$ and a constant
$C=C(\Gamma,\Lambda)$ such that
\[
  \card(A_n)\le Cn+1
\]
and
\[
  \norm p_E=\norm p_\Gamma\le \Lambda\norm p_{A_n},
  \qquad p\in\mathcal P_n(\C).
\]
The sets $A_n$ may be chosen as images under the inverse exterior conformal map of equally spaced
points of $\T$.
\end{theorem}

\begin{proof}
Set
\[
  q:=1-\Lambda^{-1},
  \qquad
  \tau:=\frac{q}{e+q}.
\]
Choose $0<c\le c_0$ so that $\omega_\Gamma(c)\le\tau$. Define
\[
  N_n:=\left\lceil\frac{\pi n}{c}\right\rceil,
  \qquad
  \theta_j:=\frac{2\pi j}{N_n},
  \qquad
  A_n:=\{\Psi(e^{i\theta_j}):0\le j<N_n\}.
\]
For every $\theta$ there is $j$ with circular angular distance
\[
  |\theta-\theta_j|_\T\le \frac{\pi}{N_n}\le\frac cn.
\]
Property (SC) then gives
\[
  |\Psi(e^{i\theta})-\Psi(e^{i\theta_j})|
  \le \tau\rho_n(\Psi(e^{i\theta})).
\]
Thus $A_n$ satisfies the hypothesis of Corollary~\ref{cor:rhonet}. Hence
\[
  \norm p_E\le \frac{1}{1-q}\norm p_{A_n}
  =\Lambda\norm p_{A_n}.
\]
Finally,
\[
  \card(A_n)=N_n\le \frac{\pi}{c}n+1.
\]
\end{proof}

\begin{remark}
Theorem~\ref{thm:SC} isolates the transferable statement. Quasicircularity is one robust
way to obtain (SC), but it is not built into the analytic argument. We do not attempt here
to characterize the Jordan curves satisfying (SC).
\end{remark}

\section{Quasicircles satisfy scale comparison}\label{sec:QC}

We now verify (SC) using the classical quasisymmetric control furnished by a
quasiconformal extension of the exterior conformal map of a quasidisk; see
\eqref{eq:QS} and \cite{AndrievskiiBelyiDzjadyk,AndrievskiiSurvey2002}. In the
form needed here, the argument reduces to a single three-point comparison.

\begin{lemma}\label{lem:QCsc}
Let $\Gamma$ be a $K$-quasicircle. Then $\Gamma$ satisfies
\textup{(SC)}. More precisely, one may take $c_0=1$ and
\[
  \omega_\Gamma(c)=\eta_K(c),
  \qquad 0<c\le1,
\]
where $\eta_K$ depends only on $K$.
\end{lemma}

\begin{proof}
Fix $n\ge1$ and $\theta\in\R$. Put
\[
  x=e^{i\theta},
  \qquad z=\Psi(x).
\]
Choose $u_n\in r_n\T$ such that
\[
  |z-\Psi(u_n)|=\rho_n(z).
\]
Such a point exists by compactness. Since $|u_n|=1+1/n$ and $|x|=1$,
\begin{equation}\label{eq:preimageradial}
  |x-u_n|\ge ||u_n|-|x||=\frac1n.
\end{equation}
Let $y=e^{i(\theta+h)}$ with $|h|\le c/n$. Then
\[
  |x-y|=|1-e^{ih}|\le |h|\le \frac cn.
\]
By \eqref{eq:preimageradial},
\[
  \frac{|x-y|}{|x-u_n|}\le c.
\]
All three points lie in $\{|w|\ge1\}$, where the quasiconformal extension agrees with
$\Psi$. Quasisymmetry therefore gives
\[
  |\Psi(x)-\Psi(y)|
  \le \eta_K(c)|\Psi(x)-\Psi(u_n)|
  =\eta_K(c)\rho_n(z).
\]
Since $\eta_K(c)\to0$ as $c\downarrow0$, this is (SC).
\end{proof}

Combining Lemma~\ref{lem:QCsc} with Theorem~\ref{thm:SC} yields the main result.

\begin{theorem}\label{thm:main}
Let $\Gamma\subset\C$ be a $K$-quasicircle and let $E$ be the closed Jordan domain
bounded by $\Gamma$. For every $\Lambda>1$ there are a constant
$C=C(K,\Lambda)$ and finite sets $A_n\subset\Gamma$ such that
\[
  \card(A_n)\le Cn+1
\]
and
\[
  \norm p_E=\norm p_\Gamma\le\Lambda\norm p_{A_n},
  \qquad p\in\mathcal P_n(\C).
\]
The sets $A_n$ may be chosen as images under the inverse exterior conformal map of equally spaced
points of $\T$.
\end{theorem}

\begin{proof}
By Lemma~\ref{lem:QCsc}, property \textup{(SC)} holds with
\[
  \omega_\Gamma(c)=\eta_K(c),
\]
where $\eta_K$ depends only on $K$.

In the proof of Theorem~\ref{thm:SC}, the quantities
\[
  q=1-\Lambda^{-1},
  \qquad
  \tau=\frac{q}{e+q}
\]
depend only on $\Lambda$. Hence $0<c\le 1$ may be chosen so that
\[
  \eta_K(c)\le\tau,
\]
with $c$ depending only on $K$ and $\Lambda$. The construction
therefore gives
\[
  \card(A_n)\le \frac{\pi}{c}n+1,
\]
so the constant $C$ depends only on $K$ and $\Lambda$.
\end{proof}

\begin{remark} The boundary correspondence
$e^{i\theta}\mapsto\Psi(e^{i\theta})$
pushes normalized Lebesgue measure on $\T$ to harmonic measure
of $\Omega$ viewed from infinity, which here coincides with the
equilibrium measure of $E$; see, for example, \cite{Ransford}. Thus the nodes in
Theorem~\ref{thm:main} divide the boundary into equal harmonic-measure pieces. Their
Euclidean spacing need not be comparable from one part of the curve to another.
\end{remark}

\section{The classical von Koch snowflake}\label{sec:koch}

The von Koch snowflake boundary $\Gamma_{\mathrm{Koch}}$ is a quasicircle.
Indeed, its standard parametrization is injective
\cite[discussion following Proposition~3.2]{vanGoldenKombrinkSamuel},
and \cite[Theorem~3.3]{vanGoldenKombrinkSamuel} yields bounded turning
for the case $p=1/3$. The characterization \eqref{eq:BT} then gives
the claim. Therefore Theorem~\ref{thm:main} gives the following.

\begin{corollary}\label{cor:koch}
For every $\Lambda>1$ there exist finite sets $A_n\subset\Gamma_{\mathrm{Koch}}$ such that
\[
  \card(A_n)=O(n)
\]
and
\[
  \norm p_{\Gamma_{\mathrm{Koch}}}\le\Lambda\norm p_{A_n},
  \qquad p\in\mathcal P_n(\C).
\]
The nodes may be chosen equally spaced with respect to harmonic measure at infinity.
\end{corollary}

The conclusion does not contradict the fractal covering geometry of $\Gamma_{\mathrm{Koch}}$. A uniform
Euclidean net sees the Hausdorff and Minkowski dimension of the snowflake, whereas the
present mesh is uniform only in the conformal angular parameter. The corresponding
Euclidean scale is the point-dependent quantity $\rho_n(z)$.

\section{Final remarks}

The proof separates two issues which are often entangled in mesh constructions.
Bernstein--Walsh and Cauchy's formula provide a local norming scale $\rho_n(z)$ without
using smoothness. The geometric problem is to count a net at that variable scale. For
quasicircles, quasisymmetry solves the counting problem through (SC).

We have made no attempt to optimize the oversampling constant $C(K,\Lambda)$. The
argument is first-order and is designed to establish linear cardinality. On the unit circle, for example, if $N>n$, then
\[
  \norm p_{\T}
  \le
  \sec\!\left(\frac{\pi n}{2N}\right)
  \max_{\omega^N=1}|p(\omega)|,
  \qquad p\in\mathcal P_n(\C);
\]
see \cite{Dubinin2011}. This estimate is sharp when $N$
is divisible by $n$.

\section*{Disclosure on AI-assisted tools}

AI-based tools were used during the exploratory phase of this project for mathematical
discussion, bibliographic assistance and drafting support. The authors independently reviewed the mathematical proofs and the
final manuscript and take full responsibility for its contents.

\end{document}